\documentclass[12pt]{amsart}
\usepackage{cleveref}
\usepackage{amssymb}
\usepackage{setspace}
\usepackage{comment}
\usepackage{fancyhdr}
\usepackage{url}

\theoremstyle{definition}
\newtheorem{thm}{Theorem}[section]
\newtheorem{lem}[thm]{Lemma}
\newtheorem{cor}[thm]{Corollary}
\newtheorem{prop}[thm]{Proposition}

\newtheorem{rem}[thm]{Remark}

\begin{document}

\title[Regular norm and operator semi norm]{Regular norm and the operator seminorm on a non-unital complex commutative Banach Algebra}

\author{Adam Orenstein}

\address{Department of Mathematics, University at Buffalo, Buffalo, NY 14260 , USA}

\email{adamoren@buffalo.edu}


\begin{abstract}
We show that if $\mathfrak{A}$ is a commutative complex non-unital Banach Algebra with norm $\|\cdot\|$, then $\|\cdot\|$ is regular on $\mathfrak{A}$ if and only if $\|\cdot\|_{op}$ is a norm on $\mathfrak{A}\oplus \mathbb{C}$ and $\mathfrak{A}\oplus\mathbb{C}$ is a commutative complex Banach Algebra with respect to $\|\cdot\|_{op}$.
\end{abstract}

\maketitle

\section{Background}\label{regBack}


Let $\mathfrak{A}$ be complex non-unital Banach Algebra $\mathfrak{A}$ with norm $\|\cdot\|$ and $\mathfrak{A}_1^+$ be the unitization of $\mathfrak{A}$ with norm $\|\cdot\|_1$ defined by $\|(a,\lambda)\|_1 = \|a\|+|\lambda|$ for all $a\in\mathfrak{A}$ and $\lambda\in\mathbb{C}$.  More about the unitization of a non-unital Banach Algebra can be found in \cite{zhuAl}.  Let $\|\cdot\|_{op}:\mathfrak{A}\oplus\mathbb{C}\rightarrow [0,\infty)$ be defined by \begin{equation}\label{opSemi}\|(a,\lambda)\|_{op}=\sup\{\|ax+\lambda x\|, \|xa+\lambda x\|:x\in \mathfrak{A}, \|x\|\leq1\}.\end{equation}  Straightforward calculations show $\|\cdot\|_{op}$ is a seminorm on $\mathfrak{A}\oplus\mathbb{C}$.  As in \cite{GaurKov}, we call $\|\cdot\|_{op}$ the operator seminorm on $\mathfrak{A}\oplus\mathbb{C}$.  If $\|\cdot\|_{op}$ is a norm on $\mathfrak{A}\oplus\mathbb{C}$, then we let $\mathfrak{A}_{op}^+$ denote the normed algebra $\mathfrak{A}\oplus\mathbb{C}$ with addition, scalar multiplication and multiplication defined as in $\mathfrak{A}_1^+$ and with norm $\|\cdot\|_{op}$.

The norm $\|\cdot\|$ on $\mathfrak{A}$ is by definition regular if for all $a\in \mathfrak{A}$, \begin{equation}\label{regDef}\|a\|=\sup\{\|ax\|,\|xa\|:x\in \mathfrak{A}, \|x\|\leq1\}\end{equation} \cite{Tak}.  For any $a\in \mathfrak{A}$, $\sup\{\|ax\|,\|xa\|:x\in \mathfrak{A}, \|x\|\leq1\}\leq\|a\|$.  This means $\|\cdot\| \text{ is regular on }\mathfrak{A} \text{ if and only if }\|a\|\leq \sup\{\|ax\|,\|xa\|:x\in \mathfrak{A}, \|x\|\leq1\}$ for all $a\in \mathfrak{A}$.  Moreover for any $a\in \mathfrak{A}$, $\|(a,0)\|_1=\|a\|$ and $\|(a,0)\|_{op}=\sup\{\|ax\|,\|xa\|:x\in \mathfrak{A}, \|x\|\leq1\}$.  So \begin{equation}\label{regCond}\|\cdot\| \text{ is regular on }\mathfrak{A} \text{ if and only if }\|(a,0)\|_1\leq\|(a,0)\|_{op}\end{equation}for all $a\in \mathfrak{A}$.

Clearly if $\mathfrak{A}$ has unit 1 with $\|1\|=1$, then $\|\cdot\|$ is regular on $\mathfrak{A}$.  But if $\mathfrak{A}$ is non-unital, then $\|\cdot\|$ may not be regular.  For consider $\ell^1$ with componentwise multiplication and let $\|\cdot\|_{\ell^1}$ be the $\ell^1$ norm.  Let $x=\left\{\frac{1}{n^2}\right\}_{n=1}^\infty$.  Then for any $y=\{y_j\}_{j=1}^\infty\in\ell^1$ with $\|y\|_{\ell^1}\leq 1$, $\|yx\|_{\ell^1}=\|xy\|_{\ell^1}=\sum_{j=1}^\infty \frac{|y_j|}{j^2}\leq \sum_{j=1}^\infty |y_j|\leq1$  It follows that $\sup\{\|xy\|_{\ell^1},\|yx\|_{\ell^1}:\|y\|_{\ell^1}\leq 1\}\leq1$.  But $\|x\|_{\ell^1}=\frac{\pi^2}{6}$.

The notion of a regular norm and its relation with $\|\cdot\|_{op}$ has been studied in \cite{Arh,GaurKov,Tak}.  In \cite{GaurKov} the following question was asked: "is the regularity of $\|\cdot\|$ equivalent to $\mathfrak{A}_{op}^+$ being a Banach Algebra?"  We will prove here that the answer is yes if $\mathfrak{A}$ is commutative.  More specifically, we will prove the following theorem.

\begin{thm}\label{mainThmReg}
Let $\mathfrak{A}$ be a complex commutative Banach Algebra with no unit and with norm $\|\cdot\|$.  Let $\mathfrak{A}_1^+$ and $\mathfrak{A}_{op}^+$ be as above.  Then $\|\cdot\|$ is regular on $\mathfrak{A}$ if and only if $\mathfrak{A}_{op}^+$ is a complex Banach Algebra with respect to $\|\cdot\|_{op}$.
\end{thm}

\section{Numerical Range}

In order to prove \Cref{mainThmReg}, we will need the notion of the numerical range of an element in a complex unital Banach Algebra.  Let $\mathfrak{B}$ be any unital complex Banach Algebra with unit $1_\mathfrak{B}$, norm $\|\cdot\|_\mathfrak{B}$ and dual space $\mathfrak{B}^{'}$.  Let $S(\mathfrak{B})=\{y\in\mathfrak{B}:\|y\|_\mathfrak{B}=1\}$.  For any $y\in S(\mathfrak{B})$ let $D(\mathfrak{B},y)=\{f\in \mathfrak{B}^{'}:\|f\|=f(y)=1\}$.  The elements of $D(\mathfrak{B},1_\mathfrak{B})$ are called the normalized states (on $\mathfrak{B}$) \cite{Bon}.  For any $b\in \mathfrak{B}$ and $y\in S(\mathfrak{B})$, the sets $V(\mathfrak{B},b)$ and $V(\mathfrak{B},b,y)$ are defined by \begin{equation}\label{numRangeDef} \begin{split}&V(\mathfrak{B},b,y)=\{f(by):f\in D(\mathfrak{B},y)\}\\& \text{ and } V(\mathfrak{B}, b)=\bigcup_{y\in S(\mathfrak{B})}V(\mathfrak{B},b,y).\end{split}\end{equation}  $V(\mathfrak{B},b)$ is called the numerical range of $b$.  Also let $\sigma_\mathfrak{B}(b)$ be the spectrum of $b\in\mathfrak{B}$ and let $\text{co}(\sigma_\mathfrak{B}(b))$ be the convex hull of $\sigma_\mathfrak{B}(b)$.

We will need the following results in order to prove \Cref{mainThmReg}.  The first two are proved in \cite{Bon}.

\begin{lem}\label{numRangProp}
Let $\mathfrak{B}$ be a unital complex Banach Algebra with unit $1_\mathfrak{B}$.  Then for any $b\in\mathfrak{B}$, $V(\mathfrak{B},b)=V(\mathfrak{B},b,1_\mathfrak{B})$.
\end{lem}

\begin{thm}\label{specNum}
Let $\mathfrak{B}$ be a unital complex Banach Algebra with unit $1_\mathfrak{B}$ and norm $\|\cdot\|_\mathfrak{B}$.  Let $\mathcal{N}_\mathfrak{B}$ be the set of all algebra norms $p$ on $\mathfrak{B}$ equivalent to $\|\cdot\|_\mathfrak{B}$ such that $p(1_\mathfrak{B})=1$ and $p(bc)\leq p(c)p(b)$ for all $c,b\in\mathfrak{B}$.  For each $p\in\mathcal{N}_\mathfrak{B}$, let $V_p(\mathfrak{B},b)$ be the numerical range of $b\in \mathfrak{B}$ with $p$ replacing $\|\cdot\|_\mathfrak{B}$.  Then for all $b\in \mathfrak{B}$ \[\text{co}(\sigma_\mathfrak{B}(b))=\bigcap_{p\in\mathcal{N}_\mathfrak{B}}V_p(\mathfrak{B},b)\]
\end{thm}

The next theorem is proved in \cite{Gol}.

\begin{thm}\label{normStates}
Let $\mathfrak{B}$ be a unital complex commutative Banach Algebra and let $f\in \mathfrak{B}^{'}$.  Then $f\in D(\mathfrak{B},1_\mathfrak{B}) \text{ if and only if } f(b)\in \text{co}(\sigma_\mathfrak{B}(b))$ for all $b\in \mathfrak{B}$.
\end{thm}

The next lemma is easy to prove.  \begin{lem}\label{numDisk}
If $\mathfrak{B}$ is a non-unital complex Banach Algebra with norm $\|\cdot\|_\mathfrak{B}$, then for all $b\in \mathfrak{B}$ and $\lambda\in\mathbb{C}$, \[V(\mathfrak{B}_1^+, (b,\lambda))=\|b\|_\mathfrak{B}\overline{\mathbb{D}}\times \{\lambda\}\] where $\|b\|_\mathfrak{B}\overline{\mathbb{D}}\times \{\lambda\}=\{(\|b\|_\mathfrak{B}w,\lambda):w\in\overline{\mathbb{D}}\}$.
\end{lem}

\section{Proof of \Cref{mainThmReg}}

\begin{proof}
$(\Leftarrow)$ Assume $\|\cdot\|_{op}$ is a norm and $\mathfrak{A}_{op}^+$ is a complex Banach algebra with respect to $\|\cdot\|_{op}$.  Note that $\mathfrak{A}_{op}^+$ is commutative as $\mathfrak{A}$ is.  Let $\Psi:\mathfrak{A}_1^+\rightarrow \mathfrak{A}_{op}^+$ be defined by $\Psi((a,\lambda))=(a,\lambda)$.  Clearly $\Psi$ is bijective, linear and since $\|(a,\lambda)\|_{op}\leq \|(a,\lambda)\|_1$ for all $(a,\lambda)\in \mathfrak{A}\oplus\mathbb{C}$, $\|\Psi\|\leq 1$.  Then by the Open Mapping Theorem \cite{Rud}, $\|\Psi^{-1}\|<\infty$.  Hence there exists $\delta>0$ so that for all $(a,\lambda)\in \mathfrak{A}\oplus\mathbb{C}$, \begin{equation}\label{normIneq} \|(a,\lambda)\|_{op}\leq \|(a,\lambda)\|_1\leq \delta\|(a,\lambda)\|_{op}.\end{equation}  Thus $\|\cdot\|_{op}$ and $\|\cdot\|_1$ are equivalent.  

Let $\mathcal{N}$ be the set of all algebra norms $p$ on $\mathfrak{A}\oplus\mathbb{C}$ equivalent to $\|\cdot\|_1$ such that $p((0,1))=1$, $p((a,\lambda)(b,\gamma))\leq p((a,\lambda))p((b,\gamma))$ for all $a,b\in\mathfrak{A}$ and $\lambda,\gamma\in\mathbb{C}$.  By \eqref{normIneq} and the fact that $\|(0,1)\|_{op}=1$, $\|\cdot\|_{op}\in \mathcal{N}$.  Then by \Cref{specNum}, \begin{equation}\label{conHullInter}\text{co}(\sigma_{{\mathfrak{A}_1}^+}(a,\lambda))\subseteq V(\mathfrak{A}_{op}^+,(a,\lambda))\end{equation} for all $(a,\lambda)\in \mathfrak{A}_{op}^+$.  Moreover by \Cref{normStates} and \Cref{numRangProp}, \begin{equation}\label{conHullSpec}V(\mathfrak{A}_1^+, (a, \lambda))= \text{co}(\sigma_{{\mathfrak{A}_1}^+}((a,\lambda)))\end{equation} for all $(a,\lambda)\in \mathfrak{A}_{op}^+$.  Thus by \eqref{conHullInter} \[ V(\mathfrak{A}_1^+, (a,\lambda))\subseteq V(\mathfrak{A}_{op}^+, (a,\lambda))\] for all $(a,\lambda)\in \mathfrak{A}_{op}^+$.  Interchanging $\mathfrak{A}_1^+$ with $\mathfrak{A}_{op}^+$ and $\|\cdot\|_1$ with $\|\cdot\|_{op}$ in the above argument yields \[V(\mathfrak{A}_{op}^+, (a, \lambda))\subseteq V(\mathfrak{A}_1^+, (a,\lambda))\] for all $(a,\lambda)\in \mathfrak{A}_{op}^+$.  Hence \begin{equation}\label{numRangeEqu} V(\mathfrak{A}_1^+, (a,\lambda))=V(\mathfrak{A}_{op}^+, (a,\lambda))\end{equation} for all $(a,\lambda)\in \mathfrak{A}_{op}^+$.

By \Cref{numDisk} \begin{equation}\label{coordZero}V(\mathfrak{A}_1^+, (a,0))=\|a\|\overline{\mathbb{D}}\end{equation} for every $a\in \mathfrak{A}$.  Also for any $F((a,0))\in V(\mathfrak{A}_{op}^+, (a,0))$, $|F(a,0)|\leq \|(a,0)\|_{op}$.  It follows from this and \eqref{coordZero} that $\|a\|\leq \|(a,0)\|_{op}$ for all $a\in \mathfrak{A}$.  That is $\|(a,0)\|_1\leq \|(a,0)\|_{op}$ for all $a\in \mathfrak{A}$.  Therefore by \eqref{regCond}, $\|\cdot\|$ is regular on $\mathfrak{A}$.

$(\Rightarrow)$ Assume $\|\cdot\|$ is regular on $\mathfrak{A}$.  As stated in \Cref{regBack}, $\|\cdot\|_{op}$ is a seminorm.  Let $(a,\lambda)\in \mathfrak{A}\oplus\mathbb{C}$ and assume $\|(a,\lambda)\|_{op}=0$.  Then for all $x\in \mathfrak{A}$ with $\|x\|\leq1$, $\|ax+\lambda x\|=0$.  Hence \begin{equation}\label{normOP}ax+\lambda x=0 \text{ and } ax=-\lambda x \end{equation} for all $x\in \mathfrak{A}$.

Suppose $\lambda\neq0$.  Then \eqref{normOP} implies $\left(\frac{a}{-\lambda}\right)x=x$ for all $x\in \mathfrak{A}$.  So $\mathfrak{A}$ has a unit.  This is a contradiction since by assumption $\mathfrak{A}$ is not unital.  Thus $\lambda=0$ and $ax=0$ for all $x\in \mathfrak{A}$.  Then since $\|\cdot\|$ is regular, $\|a\|=0$ and $a=0$.  Therefore $(a,\lambda)=(0,0)$ and $\|\cdot\|_{op}$ is a norm on $\mathfrak{A}\oplus \mathbb{C}$.

Now we will prove $\mathfrak{A}_{op}^+$ is complete.  This part is based on some calculations used in the proof of Theorem 2.3 from \cite{GaurHus}.  Let $\mathfrak{L}(\mathfrak{A})$ be the space of all bounded linear operators on $\mathfrak{A}$ with operator norm $\|\cdot\|$.  For each $(a,\lambda)\in \mathfrak{A}\oplus\mathbb{C}$, let $L_{(a,\lambda)}:\mathfrak{A}\rightarrow \mathfrak{A}$ be defined by \[L_{(a,\lambda)}(x)=ax+\lambda x.\]  It is easy to see that $L_{(a,\lambda)}$ is linear for all $(a,\lambda)\in \mathfrak{A}\oplus\mathbb{C}$ and \begin{equation}\label{mapNorm}\|L_{(a,\lambda)}\|=\|(a,\lambda)\|_{op}.\end{equation}  Thus $L_{(a,\lambda)}\in\mathfrak{L}(\mathfrak{A})$ for all $(a,\lambda)\in \mathfrak{A}_{op}^+$.

Let $\mathcal{B}=\{L_{(a,\lambda)}:(a,\lambda)\in \mathfrak{A}_{op}^+\}$ and $\Gamma:\mathcal{B}\rightarrow \mathbb{C}$ by $\Gamma(L_{(a,\lambda)})=\lambda$.  Now $\mathcal{B}$ is a subalgebra of $\mathfrak{L}(\mathfrak{A})$ and $\Gamma$ is a well-defined linear map with $\ker(\Gamma)=\{L_{(a,0)}:a\in \mathfrak{A}\}$.  Since $\|\cdot\|$ is regular on $\mathfrak{A}$, \eqref{mapNorm} tells us $\|L_{(a,0)}\|=\|a\|$.  This implies the mapping $a\mapsto L_{(a,0)}$ is an surjective isometry on $\mathfrak{A}$ and hence $\ker(\Gamma)$ is a closed subalgebra of $\mathcal{B}$.  It follows that $\Gamma$ is a bounded linear transformation.

Now let $\{(a_n, \lambda_n)\}_n$ be a Cauchy sequence in $\mathfrak{A}_{op}^+$.  Let $\epsilon>0$.  Then there exist $N_1>0$ so that \begin{equation}\label{cauchy1}n,m\geq N_1 \Rightarrow\|(a_n,\lambda_n)-(a_m,\lambda_m)\|_{op}<\epsilon.\end{equation}  Also by \eqref{mapNorm}, $\{L_{(a_n, \lambda_n)}\}_n$ is Cauchy in $\mathcal{B}$.  Thus by the above work $\{\Gamma(L_{(a_n, \lambda_n)})\}_n=\{\lambda_n\}_n$ is a Cauchy sequence in $\mathbb{C}$.  So there exist $N_2>0$ such that \begin{equation}\label{cauchy2}m,n\geq N_2 \Rightarrow|\lambda_n-\lambda_m|<\epsilon.\end{equation}  Now for any $x\in \mathfrak{A}$ with $\|x\|\leq1$, \begin{equation}\label{cauchy3}\begin{split}\|a_nx-a_mx\|&\leq \|a_nx+\lambda_nx-(a_mx+\lambda_m x)\|+\|\lambda_n x-\lambda_m x\|\\&\leq \|a_nx+\lambda_nx-(a_mx+\lambda_m x)\|+|\lambda_n -\lambda_m|.\end{split}\end{equation}  It follows that \begin{equation}\label{cauchy4}n,m\geq \max\{N_1, N_2\} \Rightarrow \|a_nx-a_mx\|<2\epsilon\end{equation} for all $x\in \mathfrak{A}$ with $\|x\|\leq1$.  Thus since $\|\cdot\|$ is regular, $\{a_n\}_n$ is Cauchy in $\mathfrak{A}$.  So $\{a_n\}_n$ converges in $\mathfrak{A}$ to some $a\in \mathfrak{A}$.  Moreover $\{\lambda_n\}_n$ converges to $\lambda\in\mathbb{C}$ as $\{\lambda_n\}_n$ is Cauchy.  Hence $\lim_{n\rightarrow\infty}\|(a_n,\lambda_n)-(a,\lambda)\|_1=0$. Therefore since $\|\cdot\|_{op}\leq \|\cdot\|_1$, $\{(a_n, \lambda_n)\}_n$ converges to $(a,\lambda)$ with respect to $\|\cdot\|_{op}$ and $\mathfrak{A}_{op}^+$ is a Banach Algebra.

\end{proof}

\section{Corollary}

The following proposition is proved in \cite{GaurKov}.

\begin{prop}
Among all the unital norms on $\mathfrak{A}\oplus\mathbb{C}$ which extend the norm $\|\cdot\|$ on $\mathfrak{A}$, $\|\cdot\|_1$ is maximal and if $\|\cdot\|$ is regular, then $\|\cdot\|_{op}$ is minimal.
\end{prop}

Combining this proposition with \Cref{mainThmReg} and \eqref{normIneq} yields the following Corollary.

\begin{cor}
If $\mathfrak{A}$ is a commutative complex non-unital Banach Algebra, then all unital norms on $\mathfrak{A}\oplus\mathbb{C}$ are equivalent $\Leftrightarrow \mathfrak{A}_{op}^+$ is a Banach Algebra $\Leftrightarrow \|\cdot\|$ is regular on $\mathfrak{A}$.
\end{cor}

\end{document}